\def\sym{\fam\comfam\com}     
\font\tensym=msbm10
\font\sevensym=msbm7 
\font\fivesym=msbm5

\newfam\symfam
\textfont\symfam=\tensym
\scriptfont\symfam=\sevensym
\scriptscriptfont\symfam=\fivesym
\def\sym{\fam\symfam\relax}

\def\R{{\sym R}}

\documentclass[12pt]{article}
\usepackage{amsfonts}
\usepackage{amsmath}
\usepackage{amssymb}
\usepackage{amsthm}
\usepackage{graphicx}
\usepackage{latexsym}
\makeatletter
\renewcommand{\@seccntformat}[1]
{\csname the#1\endcsname.\enspace}
\makeatother
\usepackage{colortbl}
\definecolor{text1}{cmyk}{1,.35,0,0} % blue text colour
\definecolor{text2}{rgb}{1,0,0} % red text colour
\definecolor{text3}{cmyk}{0,0,0,1} % black text colour
\definecolor{text4}{cmyk}{0.3,0.3,0.8,0} % grey text colour
\definecolor{text5}{cmyk}{1.0,0.0,1.0,0} % green text colour

\newcommand{\tri}{\bigtriangleup}

\newcommand{\ta}{\theta}
\newcommand{\de}{\delta}

\newtheorem{theorem}{Theorem}[section]

\newtheorem{lem}{Lemma}[section]

\newtheorem{rema}{Remark}[section]

\begin{document}

\title{On shrinkage estimation of a spherically symmetric distribution for balanced loss functions \footnote{\today}}

\author{Lahoucine Hobbad\footnote{\'Ecole Nationale des Sciences Appliqu\'ees-Marrakech, Universit\'e Cadi Ayyad, Morocco (e-mail: la.hobbad@gmail.com)},  \'{E}ric Marchand\footnote{Universit\'e de Sherbrooke, D\'epartement de math\'ematiques, Sherbrooke Qc, Canada
(e-mail: eric.marchand@usherbrooke.ca)}, and Idir Ouassou\footnote{\'Ecole Nationale des Sciences Appliqu\'ees-Marrakech, Universit\'e Cadi Ayyad, Morocco (e-mail: i.ouassou@uca.ac.ma)}}
\date{}	
	\maketitle
\vspace*{-1cm}		
\begin{abstract} 
\noindent We consider the problem of estimating the mean vector $\theta$ of a $d$-dimensional spherically symmetric distributed $X$ based on balanced loss functions of the forms: {\bf (i)} $\omega \rho(\|\de-\de_{0}\|^{2}) +(1-\omega)\rho(\|\de - \theta\|^{2})$  and {\bf (ii)} $\ell\left(\omega \|\de - \de_{0}\|^{2} +(1-\omega)\|\de - \theta\|^{2}\right)$, where $\delta_0$ is a target estimator, and where   $\rho$ and $\ell$ are increasing and concave functions.   
For $d\geq 4$ and the target estimator $\delta_0(X)=X$, we provide Baranchik-type estimators that dominate $\delta_0(X)=X$ and are minimax.  The findings represent extensions of those of Marchand \& Strawderman (\cite{ms2020}) in two directions:  {\bf (a)} from scale mixture of normals to the spherical class of distributions with Lebesgue densities  and {\bf (b)} from completely monotone to concave $\rho'$ and $\ell'$. 
\end{abstract}

\section{Introduction}

The balanced loss function (BLF) was introduced and formulated by Zellner (1994) in order  to reflect two criteria,  namely goodness of fit and precision  of estimation.   For estimating $\gamma(\theta) \in \mathbb{R}^d$ based on $X \sim f_{\theta}$, consider loss incurred by estimate $\delta$
\begin{equation}
\label{loss2}
%L_{\omega,\delta_0}(\theta,\delta)=
\omega \, \rho(||\de - \de_{0}||^{2}) \,  + \, (1-\omega) \, \rho(||\de - \gamma(\ta)||^{2})\,,
\end{equation}
where $\delta_0(X)$ is a target estimator of $\gamma(\theta)$, $\rho(\cdot) \geq 0$, and $\omega \in [0,1)$.    Zellner's original BLF corresponds to $\rho(t)=t$ and $\delta_0$ as a least-squares estimator in a regression framework.    The above loss encapsulates a more general choice of the target estimator (e.g., \cite{jjmp2006}) and a more general choice of $\rho$ (e.g., \cite{jjlm2014},  \cite{jjmp2012}), with the first term measuring proximity of estimate $\delta$ to the target $\delta_0$ in comparison to the second term measuring proximity of $\delta$ to the estimand $\gamma(\theta)$, weighted by $\omega$ and $1-\omega$ respectively.  Decision making under such a loss will necessarily lead to a compromise modulated by the amplitude of $\omega$, with the case $\omega=0$ corresponding to the so-called unbalanced loss (denoted $L_0$) where decisions are not influenced by $\delta_0$.  Alternatively, one may view the first term as a penalty for estimate $\delta$ diverging from the target $\delta_0$. An example of this arises with the choice $\delta_0=0$ connecting the balanced loss with a ridge regression or Tikhonov regularization framework. For such reasons, balanced loss functions are appealing for decision making and have interested researchers over the years.

\noindent  A natural and interesting modification of Zellner's original balanced loss function, introduced in \cite{ms2020} is given by:
\begin{equation}
\label{loss3}
%L_{\omega,\delta_0}(\theta,\delta) \, = 
\, \ell \left(\omega \, \|\de - \de_{0}||^{2} +(1-\omega) \, \|\de - \gamma(\ta)\|^{2} \right)\,,
\end{equation}
\noindent where $\delta_0$ is a target estimator of $\gamma(\theta)\in\R^d$,  $0 \leq \omega < 1 $, and $l(\cdot) \geq 0$.  Such losses possess similar attractive features as those in (\ref{loss2}), and are also especially appealing if $\ell$ is concave or even bounded.   

For the original squared error loss balanced loss function case with $\rho(t)=t$ or $\ell(t)=t$, it is known (e.g., \cite{dgs1999}, \cite{jjmp2006}, \cite{ms2020}) that frequentist risk performance of estimators is directly related to the frequentist risk performance of associated estimators under unbalanced loss.  For instance, we have the following.

\begin{lem} (Corollary 1 of \cite{ms2020})
\label{connections}
Let $X \sim f_{\theta}$ and consider the problem of estimating $\gamma(\theta)$.  Then, for $\omega \in (0,1)$,  $\delta_0(X) \, + (1-\omega) \, g_1(X)$ dominates $\delta_0(X) \, + (1-\omega) \, g_2(X)$ under loss (\ref{loss2}) with $\rho(t)=t$ if and only if $\delta_0(X) \, + \, g_1(X)$ dominates $\delta_0(X) \, + \, g_2(X)$ under loss $\|\delta-\gamma(\theta)\|^2$.
\end{lem} 

\noindent  Marchand and  Strawderman \cite{ms2020} considered the estimation of the mean of a multivariate normal $\theta$ or a scale mixture of normal distribution under losses (\ref{loss2}) and (\ref{loss3}).  They provided, for three dimensions or more, for increasing and concave $\rho$ and $l$  which also satisfy a completely monotone property, Baranchik-type estimators of $\theta$  which dominate the benchmark $\delta_0(X) = X$.  Their findings apply to a vast collections of $\rho$'s and $\ell$'s, and quite generally for scale mixtures of normal distributions, but they do not cover non completely monotone $\rho'$ and $\ell'$, as well as other spherically symmetric distributions.

\noindent  In this paper, we provide extensions with respect to the choices of $\rho$ and $\ell$, as well to the class of spherically symmetric densities.
More precisely, for $X \in \mathbb{R}^d$, 
$d \geq 4$,  with spherically symmetric Lebesgue density $ f(\| x -\theta \|^2)$, we obtain, for estimating $\theta$ for losses of types (\ref{loss2}) and (\ref{loss3}) with target $\delta_0(X)=X$, Baranchik-type estimators (\cite{bara1970}) that dominate the benchmark estimator $X$.
Our results apply to increasing and concave $\rho$ and $\ell$ and do not require complete monotonicity.  The Baranchik-type estimators studied are of the form:
\begin{equation}
\label{Baranchik}
\delta_{a,s}(X)=\left(1-a(1-\omega) \frac{s(\| X\|^2)}{\| X\|^2}\right)X \,,
\end{equation}
with $s$ twice differentiable a.e., 
\begin{equation}
\label{conditionsBaranchik}
a>0, \, 0 \leq s(\cdot) \leq 1\,,\, s(\cdot) \neq 0, s'(\cdot) \geq 0, \hbox{ and } s''(\cdot) \leq 0\,.
\end{equation}
Such estimators include simple choices like $s(t)= \frac{t}{t+b}$ with $b \geq 0$, including James-Stein estimators for $b=0$.   It is worthwhile noting that the findings here can be applied to cases where a sample $X^{(1)}, \ldots, X^{(n)}$ is drawn from  $ f(\| x -\theta \|^2)$ and when inference is then based on an  estimator (such as $\bar{X}$) which is a spherically symmetric and translation invariant function of these data (e.g., \cite{bs1980}).
  
\noindent The paper is organized as follows.
In Sections 2 and 3, we present dominance results for balanced losses (\ref{loss2}) and (\ref{loss3}) respectively.  Both sections include a subsection of illustrations and remarks.   Several technical results are presented in these sections, and others are relegated to an Appendix. Finally, a detailed example is presented in Section 4, with a comparison of frequentist risks and some details worth sharing about the calculations themselves.

\section{ Risk analysis for loss $\omega\rho(\|\delta-X\|^{2}) +(1-\omega) \rho(\|\delta-\theta\|^{2})$}
\label{sect2}

\subsection{Preliminary results and definitions}
\label{subsection2.1}	
We begin with an initial observation. It follows from \cite{kiefer1957} that $\delta_0(X)=X$ is minimax for estimating $\theta$, whenever $X$ has density $f(\|x-\theta\|^2)$, $\theta \in \mathbb{R}^d$, and loss $\rho(\| \delta-\theta\|^{2})$ with $\rho(0)=0$ and $\rho'(t)>0$ for $t>0$, as long as $X$ has finite risk.   Furthermore, as put forth in \cite{jjmp2012},
if an estimator $\delta_{0}(X)$ is minimax for estimating $\theta$ under loss  $\rho(\| \delta-\theta\|^{2})$ then, it is also minimax under  balanced loss in (\ref{loss2}) with $\gamma(\theta)=\theta$ for all $0 < \omega <1$.  Therefore, the estimators of this section which dominate the benchmark $\delta_0(X)=X$ under balanced loss (\ref{loss2}) are minimax in cases of spherically symmetric density models and for the conditions {\bf C1} below in (\ref{C1}) on $\rho$.	
%\begin{theorem} 
%\label{t1} 
%\begin{enumerate}
%\item[{\bf (a)}]
%Let  $X \sim f_{\theta}$. If  the estimator $\delta_{0}(X)$ is minimax  under loss  $\rho(\| \delta-\theta\|^{2})$ then, it is also minimax under  balanced loss in (\ref{loss2}) with $\gamma(\theta)=\theta$ for all $0 < \omega <1$. 
%\item[{\bf (b)}]   Suppose $X \sim f(\|x-\theta\|^2)$ with non-increasing $f$.   %Then, assuming its risk is finite, $\delta_0(X)=X$ is a minimax estimator of $%\theta$ under loss $\rho(\| \delta-\theta\|^{2})$
%with $\rho'(t) \geq 0$ for all $t>0$ and $\rho'(t) > 0$ for some interval.
%\end{enumerate}
%	\end{theorem}

\noindent For the function $\rho$ in loss (\ref{loss2}), we assume the following throughout  
\begin{equation}
\label{C1}
	{\bf C1}:   \rho (0) = 0, 0<\rho^{\prime}(0) < +\infty\,, \hbox{and }\rho \hbox{ is concave} .
\end{equation}
\noindent 		Examples for which $\rho$ satisfies condition {\bf C1}, other than $\rho(t)=t$, include: {\bf (i)} $ \, \rho(t)= 1-\exp(- \alpha \, t)$ with $\alpha > 0$, {\bf (ii)} $ \,
\rho(t)=\log(1+t)$,  {\bf (iii)} $ \; \rho(t)=(1+t/\gamma)^{\beta}$ with $\gamma >0, \beta \in (0,1)$, {\bf (iv)} $\; \rho(t) = r^{2}t/(rt+1)$ with $r >0$, {\bf (v)} $\rho(t) \, = \, \arctan(t)$, and {\bf (vi)} $\rho(t) \, = \, \tanh (t)$.   Except for {\bf (v)} and {\bf (vi)}, the examples were presented by Marchand \& Strawderman in \cite{ms2020} as examples of completely monotone $\rho'$ (i.e., $(-1)^{k+1} \, \rho^{(k)}(t) \geq 0$ for $t>0$ and $k \in \mathbb{N}$).     We do not require completely monotonicity and condition {\bf C1} is weaker.    Another class of loss functions satisfying the above conditions are given by: {\bf (vii)} $\rho(t)\,=\,G(t)$ with $G$ a 
cdf on $(0,\infty)$ with non-increasing density $G'$.  Such a class of loss functions arose recently in a predictive density estimation framework described in \cite{kms2017}.  Finally, we point out that condition {\bf C1} implies:
\begin{center}
	{\bf C2}:   $\rho (0) = 0$, $0<\rho^{\prime}(0) < +\infty$, $\rho$ is concave, and $\rho(t)/t$ is non-increasing for $t \in \mathbb{R}_+$.
\end{center}
The added non-increasing property of $\rho(t)/t$ is thus superfluous, but we will use this on several occasions, as well as the outright concavity of $\rho$.

\noindent We make use throughout of standard definitions and properties of spherically symmetric distributions, as well as shrinkage estimation techniques and properties of superharmonic functions.    Some key features are either collected in this subsection or in the Appendix, and further properties and definitions can be found for instance in  \cite{fsw2018}.   

\noindent   A well-known, useful property and characterization of spherically symmetric distributions (e.g., Theorem 4.1 in \cite{fsw2018}) is the independence between the radius 
$R=\|X-\theta\|$ and $\frac{X-\theta}{\|X-\theta\|}$, with the conditional distribution of $X|R=r$ uniformly distributed (noted $U_{r,\theta}$) on the sphere $S_{r,\theta}= \{ x\in\R^d \mid \| x-\theta\|= r\}$ of radius $r$ centered at $\theta$.   Our findings capitalize on a corresponding conditional on $R$ risk decomposition.  Moreover, 
the distribution of $R$ is independent of $\theta$ and its density (called radial) is given by $h(r)\,=\, \frac{2 \pi^{d/2}}{\Gamma(d/2)} \, r^{d-1} \, f(r^2) \, \mathbb{I}_{\mathbb{R}_+}(r)$ when $X$ has Lebesgue density $f(\|x-\theta\|^2)$.
%\begin{lem} 
%\label{teo-ss-rayon} (See \cite{fsw2018}, page 128.)
%A distribution $P$ in $\R^d$ is spherically symmetric about $\theta\in\R^d$ if and only if there exists a distribution $\gamma$ in $\R_+$ such that $P(A)=\int_{\R_+}U_{r,\theta} (A) d\gamma(r)$ for any Borel set $A$ of $\R^d$. Furthermore, if a random vector $X$ has such a distribution $P$, then the radius $\|X-\theta\|$ has distribution $\rho$ (called the radial distribution) and
%the conditional distribution of $X$ given $\|X-\theta\| = r$ is the uniform %distribution $U_{r,\theta}$ on the sphere $S_{r,\theta}$ of radius $r$ and %centered at $\theta$.
%\end{lem} 

\noindent 	  Properties of superharmonic functions also play an important role. We recall that a continuous function $g:\R^{d}\longrightarrow \R $ is superharmonic if and only if: for all $t_{0} \in \mathbb{R}^d$ and $r > 0$, the average of $g$ over the surface of the sphere $S_{r, t_0}$ is less or equal than $g(t_{0})$. For twice-differentiable $g$, the superharmonicity of $g$ is equivalent to its Laplacian being less or equal to $0$, i.e., $\Delta(g)\leq 0$  with  $\Delta(g) =\sum\limits_{k=1}^{d}\left({\partial^{2} g(t)}/{\partial t_{k}^{2}}\right)$. 
As stated in Lemma \ref{lemA.3}, superharmonicity of  $g$ can be used to obtain an inequality relating the conditional expected value on the ball to the conditional expected value on the sphere. An important fact here is that if $g$ is superharmonic, its average over the ball (``volume") is greater than its average over the sphere (``surface area").   We conclude this subsection with a pivotal inequality, which appeared in \cite{ms2020} for scale mixtures of normals, but which holds here more generally over the class of spherically symmetric densities.

\begin{lem} 
\label{lm3-1}
Suppose that $X$ is spherically symmetric distributed about $\theta$ with density $f(\|x-\theta\|^2)$, that $\mathbb{E}_{\theta} (\|X\|^{-2}) < \infty$, and that the function $\rho$ satisfies \textbf{C1}.
		For $d \geq 4$, and for $s:\R^{d}\longrightarrow [0,1]$  a twice-differentiable, non-decreasing and concave function then  
		\begin{equation}
		\label{eq-lm3-1}
		E_{\theta}\left[ \rho\left(\frac{s(\| X\|^2)}{\| X\|^2}\right)\right]\leq \rho^{\prime}(0) E_{\theta}\left[ \frac{s(\| X\|^2)}{\| X\|^2}\right]
		\leq \rho^{\prime}(0)E_{\theta}\left[\frac{s(\| Y\|^2)}{\| Y\|^2}\right] \,,
		\end{equation}
%\noindent 	where $ Y\sim f^{\ast}(||y-\theta||^{2})\,=\, \frac{\rho'(||y-\theta||^{2})) \,
%f(||y-\theta||^{2})}{K} $ with $K\,=\, \mathbb{E}_0(\rho'(\|X\|^2))\,$.
%$K =\int \rho^{\prime}(r^{2})f(r^{2})dr$ and $f^{\star}(r^2)=\frac{\rho^{\prime}%(r^{2})f(r^{2})}{K}$ .
where $ Y\sim f^{\ast}(||y-\theta||^{2})$ with $ f^{\ast}(t)\,=\, \frac{\rho'(t) \, f(t)}{K}$ and $K\,=\, \mathbb{E}_0(\rho'(\|X\|^2))\,$.

%$\,=\, \frac{\rho'(||y-\theta||^{2})) \,
%f(||y-\theta||^{2})}{K} $ with $K\,=\, \mathbb{E}_0(\rho'(\|X\|^2))\,$.
%$K =\int \rho^{\prime}(r^{2})f(r^{2})dr$ and $f^{\star}(r^2)=\frac{\rho^{\prime}%%(r^{2})f(r^{2})}{K}$ . 
		%\end{enumerate}
	\end{lem}
	
\noindent {\bf{Proof.}}
The first inequality of (\ref{eq-lm3-1}) follows, on taking expectations, from the concave function inequality
	$$\rho(t)\leq \rho(0)+\rho^{\prime}(0)t=\rho^{\prime}(0)t$$
	given that $\rho$ is concave with $\rho(0)=0$.
	Denoting by $h$ and $h^{\ast}$ the densities of $\|X-\ta\|$ and $\|Y-\ta\|$ respectively, $\mathbb{E}^*$ the expectation with respect to $h^*$, and making use of the equality of the conditional distributions $X|\|X-\theta\|=r$  and $Y|\|Y-\theta\|=r$, we have 
	\begin{eqnarray}
\nonumber		E_{\theta}\left[\frac{s(\| X\|^2)}{\| X\|^2}\right]
%		&=& E\left[E_{R,\theta}\left(\frac{s(\|X\|^{2})}{\|X\|^{2}}\; %\Big{\arrowvert} \;\|X-\theta\|=R\right)\right]\\
\nonumber		&=& \int_{\R_{+}}E_{\theta}\left( \frac{s(\|X\|^{2})}{\|X\|^{2}}\;\Big{\arrowvert}\; \|X-\theta\|=r\right)h(r)dr \\
%		&=&\int_{\R_{+}}E_{r,\theta}\left( \frac{s(\|X\|^{2})}{\|X\|^{2}}\;\Big{\arrowvert}%\; \|X-\theta\|=r\right) \frac{2\pi^{d/2}}{\Gamma\left(\frac{d}{2}\right)}r^{d-1}%f(r^{2})dr \\
%		&=&\int_{\R_{+}}E_{r,\theta}\left( \frac{s(\|X\|^{2})}{||X||^{2}}\;\Big{\arrowvert}%\; ||X-\theta||=r\right)		\frac{2\pi^{d/2}}
%		{\Gamma\left(\frac{d}{2}\right)} r^{d-1}\frac{\rho^{\prime}(r^{2})f(r^{2})}{K}%\frac{K}{\rho^{\prime}(r^{2})}dr\\
\nonumber		&=&\int_{\R_{+}}E_{\theta}\left( \frac{s(\|Y\|^{2})}{\|Y\|^{2}}\;\Big{\arrowvert}\; \|Y-\theta\|=r\right)\frac{K}{\rho^{\prime}(r^{2})} \, h^{\ast}(r)dr \\
\label{e}		&=& E^{\star}\left[ E_{\theta}\left( \frac{s(\|Y\|^{2})}{\|Y\|^{2}}\;\Big{\arrowvert}\; \|Y-\theta\|=R\right)\frac{K}{\rho^{\prime}(R^{2})}\right] \,.
	\end{eqnarray}
	 % $$K =\int \rho^{\prime}(R^{2})f(R^{2})dR$$
	%$$ Y\sim f^{\ast}(||y-\theta||^{2})$$and  
%	$h^{\ast}(r)=\left({2\pi^{d/2}}	/{\Gamma\left(\frac{d}{2}\right)}\right) r^{d-1}%f^{\star}(r^2)$ is a  density function of $\|Y-\ta\|$.\\
	%$$f^{\star}(R^2)=\frac{\rho^{\prime}(R^{2})f(R^{2})}{K}$$
Since the function $s(\|t\|^{2})/\|t\|^{2}$ is superharmonic for $d \geq 4$ as shown in Lemma \ref{lem A.1}, and since for such superharmonic functions the sphere mean is non-increasing in the radius (e.g., \cite{fsw2018}, Theorem A.4, page 304), it follows that 
	$E_{\theta}\left( \frac{s(\|Y\|^{2})}{||Y||^{2}}\;\Big{\arrowvert}\; \|Y-\theta\|=r\right)$   is 
	non-increasing in $r$. Finally, along with the concavity of $\rho$ which implies that ${K}/{\rho^{\prime}(r^{2})}$ is non-decreasing in $r$, the covariance inequality applied to (\ref{e}) implies that  
	
	\begin{eqnarray*}
		 	E_{\theta}\left[\frac{s(\| X\|^2)}{\| X\|^2}\right] \, \leq \, 	E_{R}^{\star}\left[ \frac{K}{\rho^{\prime}(R^{2})}\right]
		E^{\star}\left[E_{\theta}\left( \frac{s(\|Y\|^{2})}{\|Y\|^{2}}\;\Big{\arrowvert}\; \|Y-\theta\|=R\right)\right]
		%&=&E_{r}^{\star}\left[E_{r,\theta}\left( \frac{r(\|Y\|^{2})}{\|Y\|^{2}}\;\Big{\arrowvert}\; \|Y-\theta\|=R\right)\right]
		= E_{\theta} \left[\frac{s\left(\| Y\|^2\right)}{\| Y\|^2}\right],
	\end{eqnarray*}
due to the fact that $	E_{R}^{\star}\left[ {K}/{\rho^{\prime}(R^{2})}\right]=1$.  \qed

	 \subsection{Dominance finding}
	\label{sect2.2}
We are now ready for the main dominance finding of this section.	
	\begin{theorem}
		\label{theorem-rho-lebesgue}
		Suppose that $X$ is spherically symmetric distributed about $\theta$ with density $f(\|x-\theta\|^2)$, that both $\mathbb{E}_0(\|X\|^2)$ and  $\mathbb{E}_0(\|X\|^{-2})$ are finite, and that the function $\rho$ satisfies \textbf{C1}.
		For $d\geq 4$ and for estimating $\theta$ under loss (\ref{loss2}) with $\delta_0(X)=X$, the estimator $\delta_{a,s}(X)$  in (\ref{Baranchik}) satisfying conditions (\ref{conditionsBaranchik})  dominates $\delta_0$ provided:
\begin{equation} 
\label{cutoff1}
0<a  < \frac{2K(d-2)/d}{\left\lbrace \, \omega\rho^{\prime}(0)+K(1-\omega)\right\rbrace \, E^{\star}_{0}\left(\|Y\|^{-2}\,\right)}\,,
\end{equation} 
where $E^{\star}_{\theta}$ is the expectation taken with respect to $ Y\sim f^{\ast}(||y-\theta||^{2})$ with $ f^{\ast}(t)\,=\, \frac{\rho'(t) \, f(t)}{K}$ and $K\,=\, \mathbb{E}_0(\rho'(\|X\|^2))\,$.	An equivalent condition for the above dominance condition is:
\begin{equation} 
\label{cutoff2}
0<a  < \frac{2\, K^2 \, (d-2)/d}{\left\lbrace \, \omega\rho^{\prime}(0)+K(1-\omega)\right\rbrace \, E_{0}\left(\frac{\rho'(\|X\|^{2})}{\|X\|^{2}}\,\right)}\,\,.
\end{equation}
\end{theorem}
\noindent {\bf{Proof.}} Set $\delta_g(X) \, = \, X + (1-\omega) g(X)$ and consider the  difference in risks  
$$\Delta{\cal{R}}(\ta)=
{\cal{R}}(\theta, \delta_{g})-{\cal{R}}(\theta, \delta_{0})\,, \, \theta \in \mathbb{R}^d\,. $$
We have 
	\begin{eqnarray}
	\label{dif-risk}
	\Delta{\cal{R}}(\ta)
%		&=&
%		E_{\ta}\left\lbrace\omega\rho(\|X+(1-\omega)g(X)-X\|^{2})+(1-\omega)\rho(\|X+(1-\omega)g(X)-\ta\|^{2})\right.\nonumber\\
%		&&\left.-(1-\omega)\rho(\|X-\ta\|^{2})\right\rbrace\nonumber\\
		&=&
		E_{\ta}\left\lbrace\omega\rho(\|(1-\omega)g(X)\|^{2})+(1-\omega)\rho(\|X-\ta+(1-\omega)g(X)\|^{2})\right.\nonumber\\
		&&\left. 
		-(1-\omega)\rho(\|X-\ta\|^{2})\right\rbrace\nonumber\\
		&\leq &
		E_{\ta}\left[\omega\rho^{\prime}(0)\|(1-\omega)g(X)\|^{2}+(1-\omega)\left\lbrace\rho(\|X-\ta\|^{2}+2(1-\omega) \, g(X)^{\top}(X-\theta)\right.\right.\nonumber\\
		&&\left.\left.\qquad
		+(1-\omega)^{2}\|g(X)\|^2)-\rho(\|X-\ta\|^{2})\right\rbrace\right]\nonumber\\
%		&\leq&
%		E_{\ta}\left[\omega\rho^{\prime}(0)(1-\omega)^{2}\|g(X)\|^{2}+(1-\omega)\rho^{\prime}(\|X-\ta\|^{2})\left\lbrace 2(1-\omega)\right.\right.\nonumber\\
%		&&\left.\left.\qquad
%		g(X)^{\top}(X-\theta)+(1-\omega)^{2}\|g(X)\|^2\right\rbrace\right]\nonumber\\
		& \leq &
		E_{\ta}\left[\omega\rho^{\prime}(0)(1-\omega)^{2}\|g(X)\|^{2}+(1-\omega)K\times\frac{\rho^{\prime}(\|X-\ta\|^{2})}{K}
		\left\lbrace2(1-\omega)g(X)^{\top}(X-\theta)\right.\right.\nonumber\\
		&&\left.\left.+(1-\omega)^{2}\|g(X)\|^2\right\rbrace\right]\nonumber\\
		&=&
		(1-\omega)^{2}\left[ E_{\ta}\left[\omega\rho^{\prime}(0)\|g(X)\|^{2}\right]+
		KE_{\theta} \left\lbrace2g(Y)^{\top}(Y-\theta)+\
		(1-\omega)\,\|g(Y)\|^2\right\rbrace\right] \,,
	\end{eqnarray}
where the two inequalities follow from the concave function inequality:
$\rho(t_{1})-\rho(t_{2}) \leq \rho'(t_{1})(t_{1}-t_{2}),\, $
with the condition $\rho(0)=0$.   Now, for the Baranchik  $g(x)=-a\left({s(\| x\|^2)}/{\| x\|^2}\right)x$ and $a>0$,   it is easy to show that
 $$div\left(g(x)\right)\, = \, -2as'(\| x\|^2) 	-a(d-2)\displaystyle\frac{s(\| x\|^2)}{\| x\|^2} 
 \, \leq  -a(d-2)\displaystyle\frac{s(\| x\|^2)}{\| x\|^2} \,.$$
 
\noindent  From (\ref{dif-risk}), using the condition $0 \leq s(\cdot) \leq 1$ as well as   Lemma \ref{lm3-1} and  Lemma \ref{lem3.2}, we obtain 
	\begin{eqnarray}
	\label{risk2}
		\displaystyle\frac{
		\Delta{\cal{R}}(\ta)
			}{(1-\omega)^{2}}
		&\leq& E_{\ta}\left[\omega\rho^{\prime}(0)\frac{a^{2}s(\| X\|^2)}{\| X\|^2}\right]
		+KE_{\ta}\left[2g(Y)^{\top}(Y-\theta)\right] +KE_{\ta} \left[(1-\omega)\frac{a^{2}s(\| Y\|^2)}{\| Y\|^2}\right]\nonumber\\
		&\leq &
		E_{\ta}\left[\omega\rho^{\prime}(0)\frac{a^{2}s(\| Y\|^2)}{\| Y\|^2}\right] -2 Ka\displaystyle\frac{d-2}{d}\displaystyle\int_{\R_{+}} r^{2}\displaystyle\int_{B_{r,\theta}}\displaystyle\frac{s(\| y\|^2)}{\| y\|^2}dV_{r,\theta}(y) \, h^*(r) \, dr\, \\
		&&
	+ \,	K  \, E_{\ta} \left[ (1-\omega)\frac{a^{2}s(\| Y\|^2)}{\| Y\|^2}\right] \,, \nonumber\\
		&&, 	\end{eqnarray}
$h^*$ being the radial density of $\|Y-\theta\|$.
As $t \longrightarrow {s(\| t\|^2)}/{\| t\|^2}$  is superharmonic, then according to  
	Lemma \ref{lm3.3} expression  (\ref{risk2})  is bounded above by
		\begin{eqnarray*}
		&& E_{\ta}\left[\omega\rho^{\prime}(0)\frac{a^{2}s(\| Y\|^2)}{\| Y\|^2}\right] -2 Ka\displaystyle\frac{d-2}{d}\displaystyle\int_{\R_{+}} r^{2}\displaystyle\int_{S_{r,\theta}}\displaystyle\frac{s(\| y\|^2)}{\| y\|^2}dV_{r,\theta}(y) \, h^*(r) \, dr\, \\
		&&
	+ \,	K  \, E_{\ta} \left[ (1-\omega)\frac{a^{2}s(\| Y\|^2)}{\| Y\|^2}\right] \,, \nonumber \\
%		&=&  aE\left\{E_{\ta} \left[\frac{s(\| Y\|^2)}{\| Y\|^2}\Big{\arrowvert}\|Y-\ta\|=R\right] \left( a\omega\rho^{\prime}(0)-2K\displaystyle\frac{d-2}{d}R^{2}+
%		K \, (1-\omega) \, a\right) \right\}\nonumber\\
		&=& a \, E \left\{E_{\ta} \left( R^{2}\frac{s(\| Y\|^2)}{\| Y\|^2}\Big{\arrowvert}\|Y-\ta\|=R\right) \, \left( a \, \frac{\omega\rho^{\prime}(0)+K(1-\omega)}{R^{2}}-2K\displaystyle\frac{d-2}{d}\right) \right).\nonumber \\
	\end{eqnarray*}
	
	\noindent With $E_{\ta} \left[ R^{2} \frac{s(\| Y\|^2)}{\| Y\|^2} \, \Big{\arrowvert}\|Y-\theta\|=R \right]$  non-decreasing in $R$ by Lemma \ref{lemBS1978}, and  $\left\{a\left[\omega\rho^{\prime}(0)+K(1-\omega)\right]/R^{2}\right\}$ non-increasing in $R$, an application of the covariance inequality leads to 
	\begin{eqnarray*}
\frac{\Delta{\cal{R}}(\ta)} {(1-\omega)^{2}}		&\leq&	aE \left\{E_{\ta}\left[R^{2}\frac{s(\| Y\|^2)}{\| Y\|^2}\Big{\arrowvert} \|Y-\ta\|=R\right]\right\} \times E \left[\left( a\frac{\omega\rho^{\prime}(0)+K(1-\omega)}	{R^{2}}-2K\displaystyle\frac{d-2}{d}\right)\right] \\
		&=&aE_{\theta} \left( \|Y-\ta\|^{2}\frac{s(\| Y\|^2)}{\| Y\|^2} \right) \times \left[ a \left\lbrace \omega\rho^{\prime}(0)+K(1-\omega)\right\rbrace \, E_{\theta}\left(\frac{1}{\|Y-\theta\|^{2}}\right)  -2K\displaystyle \, \frac{d-2}{d}\right] \,,\nonumber
	\end{eqnarray*}
and the result follows.  	\qed

	 \subsection{Examples and discussion}

\noindent  The dominance finding of Section 2.2 is applicable to many choices of $\rho$ and the underlying spherically symmetric $f$.  The results guarantee the existence of a class of Baranchik estimators which are minimax.  For $d \geq 4$, subject to finiteness of risk or the minimaxity of $X$, the results apply quite generally with respect to the choice of the Baranchik estimator, the weight $\omega$, the choice of $\rho$, and the underlying model density $f$.   Our findings do expand on existing results, namely the findings of Marchand and Strawderman \cite{ms2020}, with both a wider class of loss functions and underlying models.  However, they do not duplicate their results when applicable as our conditions turn out to be stronger, i.e., the cut-off points are smaller (see Section 2.3.2).  Although the focus of this paper is on the balanced case with $\omega>0$, Theorem \ref{theorem-rho-lebesgue} applies for the case $\omega=0$ nevertheless with the cut-off point in (\ref{cutoff1}) simplifying to
$2 (d-2)\,  \{ dE^{\star}_{0}\left(\|Y\|^{-2}\,\right)\,\}^{-1}$, and matching the one obtained by Brandwein and Strawderman (\cite{bs1980}, Theorem 2.1). 

\noindent  Without dwelling too much on the wealth of applicable situations or calculations of the cut-off points, we elaborate here a little bit with some illustrations and remarks.

\subsubsection{Choices of $\rho$ and determination of the cut-off points}

\noindent  Cut-off point (\ref{cutoff2}) of Theorem \ref{theorem-rho-lebesgue} is explicitly represented and is conveniently expressible in terms of the expectations $K= \mathbb{E} (\rho'(W))$ and  $\mathbb{E} (\rho'(W)/W)$, with respect to the density of $W=\|X\|^2$ when $\theta=0$ and given by 
\begin{equation}
\label{pdfW}
\frac{\pi^{d/2}}{\Gamma(d/2)} \, w^{d/2-1} \, f(w) \, \mathbb{I}_{\mathbb{R}_+}(w).
\end{equation}
With $\rho(t) \to \frac{\rho(t)}{\rho'(0)}$; which tells us that we can set $\rho'(0)=1$ without loss of generality; and by setting 
\begin{equation}
\nonumber
\noindent   I_{k} \, = \, \int_0^{\infty}  w^{k - 1} \, \rho'(w) \, f(w) \, dw\,, 
\end{equation}
cut-off point (\ref{cutoff2}) can simply be expressed as 
\begin{equation}
\label{cutoff4}
a_0 \, = \,   \frac{ 2 (\frac{d-2}{d}) \, I_{d/2}^2}{  \left\lbrace\, \omega \, \Gamma(d/2)/ \pi^{d/2}  + \, (1-\omega) I_{d/2}   \right\rbrace  \, I_{d/2 -1}} \,.
\end{equation}

 Marchand \& Strawderman provide a similar cut-off point to (\ref{cutoff2}) applicable to completely monotone $\rho'$ and a scale mixture of normals model density, while the findings of Section 2.2 apply to non-completely monotone choices of $\rho'$ such as those referred at the beginning of Section 2.1.  As an example, consider the $\rho(t) =  2 \Phi(t) - 1$ with $\Phi$ the standard normal cdf (this is the cdf of a truncated standard normal  distribution).   Then, Theorem \ref{theorem-rho-lebesgue} applies and can one simply take $\rho'(t)/\rho'(0)= e^{-t^2/2}$ for evaluating (\ref{cutoff2}) or (\ref{cutoff4}) numerically.  Even in the normal case with $X \sim N_d(\theta, I_d)$ and  $W \sim \chi^2_d(0)$, the result is new and yields $a_0$ evaluated with 
\begin{equation}
\nonumber I_{k} \, = \int_0^{\infty} w^{k-1} \, e^{-(w+w^2)/2} \, dw \,, \hbox{ for } k\,=\, d/2, d/2-1 . 
\end{equation} 
  
\subsubsection{Scale mixtures of normal distributions case}	
Theorem \ref{theorem-rho-lebesgue} applies to scale mixtures of normals admitting representation
\begin{equation}
\label{mixtureofnormals}
X|V \sim N_d(\theta, V I_d)\,, V \sim G\,,
\end{equation}
where $G$ is a c.d.f. for $V$, and including many familiar distributions such as normal, logistic, Laplace, exponential power, Student, etc.    For choices of $\rho$ satisfying condition {\bf C1} such that also $\rho'$ is completely monotone, the cut-off point given in \cite{ms2020} is greater than (\ref{cutoff2}) by a factor of $d/(d-2)$.  So, Theorem \ref{theorem-rho-lebesgue}'s result is weaker.  However, as reviewed above in subsection 2.3.1, Theorem \ref{theorem-rho-lebesgue} applies for more general $\rho$.  Also, the results of Section 2.2 apply for more general models and we pursue with such examples.    

\subsubsection{Uniform distribution on a ball}	 

\noindent Theorem \ref{theorem-rho-lebesgue} applies for the case of a uniform distribution on a ball centered at $\theta$ of radius $m$: $B_{m,\theta} \, = \, \{x \in \mathbb{R}^d: \|x-\theta\| \leq m   \}$, with $d \geq 4$ and densities
\begin{equation}
\label{density-uniform}
f(\|x-\theta\|^2) \, = \,    \frac{\Gamma(d/2+1)}{m^d \, \pi^{d/2}} \,\, \mathbb{I}_{(0,m)}(\|x-\theta\|)\,.  
\end{equation}
In such cases, $\frac{W}{m^2} \, = \, \frac{\|X-\theta\|^2}{m^2}$ has Beta density $\frac{d}{2} \, t^{d/2-1} \, \mathbb{I}_{(0,1)}(t)$, and cut-off point (\ref{cutoff2}) is readily available by either analytical or numerical evaluations of $I_{k} \,=\, \frac{\Gamma(d/2+1)}{m^d \pi^{d/2}} \int_0^{m^2}   w^{k-1} \, \rho'(w) \,   dw$ for $k=d/2$ and $d/2-1$.  As an illustration, the choice $\rho(t)=\log(1+t), d=4$  yields:
$$   a_0 \, = \,  \frac{\{m^2- \log (1+m^2)   \}^2}{\{\omega m^4/2 + (1-\omega) \, (m^{2} - \, \log(1+m^2) ) \} \log(1+m^2) }$$ 
For the unbalanced squared error loss case with $\omega=0$ and $\rho(t)=t$, early shrinkage estimation analysis was provided in \cite{b+s1978} for such uniform distributions and their mixtures, the latter playing a key role in modelling since their ensemble spans the entire class of  unimodal spherically symmetric distributions.
	 
\subsubsection{Kotz  type  distribution}
	
\noindent  Theorem \ref{theorem-rho-lebesgue} applies to Kotz model densities $f(\|x-\theta\|^2)$ with
	 \begin{equation}
	 \label{model-kotz}	 
	 f(t)\,= \, c_d \; t^{s\nu - \frac{d}{2}} \, e^{-r t^s} \,\, \mathbb{I}_{\mathbb{R}_+}(t)\,,
	 \end{equation}
with $r,s,\nu >0$ and $c_d\,=\, \frac{s \, \Gamma(\frac{d}{2})}{\pi^{\frac{d}{2}} \, \Gamma(\nu)} \, r^{\nu}$.  \footnote{A more frequent parametrization in the literature has $N=s\gamma+1-d/2$.}  The distribution, originally introduced in \cite{kotz1975} for $s=1$ has generated much interested over the years (e.g. \cite{Nadarajah2013}), namely for its flexibility in representing non-unimodal densities (for $2s\nu > d$), non-scale mixture of normal densities (any choice except $0 < s \leq 1$ and $2 s \nu = d$), as well as including the normal case ($s=1, \nu=d/2$), scale mixtures of normal ($0 < s \leq 1$ and $2 s \nu = d$), and exponential power densities ($2s\nu = d$).
Moreover, it is simple to verify that the distribution of $W^s=\|X-\theta\|^{2s}$ is distributed as Gamma$(\nu,r)$, which facilitates the expression of cut-off point (\ref{cutoff2}).  Resulting integrals will not be available in general in closed forms.  One exception arises for $s=1$ and reflected normal $\rho(t)= 1 - e^{-\alpha \, t}$.   For $\nu>1$ (we need this for finiteness of risk), an evaluation of (\ref{cutoff2}) or (\ref{cutoff4}) yields:
$$   a_0 \, = \,   \frac{2 (\nu-1) \, (1-2/d)/(r+\alpha) }{ \omega \, (1+\alpha/r)^{\nu} \, + \, (1-\omega) \, } \, .$$
		 
\noindent A numerical illustration with an underlying Kotz density is provided and commented upon in Section 4.3.

\section{Risk analysis for loss $\ell \left(\omega\|\delta-X\|^{2} +(1-\omega) \|\delta-\theta\|^{2} \right)$}

\subsection{Dominance finding}	
\label{subsection3.1}
For spherically symmetric model densities, we evaluate the frequentist risk performance of an estimator $\delta(X)$  of $\theta$ under the balanced loss (\ref{loss3}) which incorporates the target estimator $\delta_{0}(X)\, = \, X$. For the function $\ell$, we work with the following conditions  throughout this section:
\begin{equation}
\label{C3}
\textbf{C3} :   \ell (0) = 0, \ell^{\prime}(\cdot) > 0 \,, \hbox{and }\ell \hbox{ is twice-differentiable and concave} .
\end{equation}
As for the choice of $\rho$ in the previous section, the completely monotone requirement on $\ell'$ in \cite{ms2020} is relaxed here.   
Examples of $\rho$ that satisfy condition {\bf C1}, including {\bf (i)} to {\bf (vii)} following (\ref{C1}),	provide examples of $\ell$ that satisfy above condition {\bf C3}.   But, we do not require finiteness of $\ell'(0)$ so that many other losses, such as $L^{q/2}$ losses with {\bf (viii)} $\ell(t)\,=\, t^q\,, 0<q<1$, satisfy {\bf C3} as well.   Another interesting choice comes from \cite{kms2017} as: {\bf (ix)} $\ell(t) \, = \, 2 Q(\sqrt{t}) \, - \, 1$ with $Q$ a c.d.f. on $\mathbb{R}$ with even and unimodal density $Q'$.      This arises as an intrinsic loss in measuring the distance between estimate $\hat{\mu}$ and location parameter $\mu$ through the $L^1$ discrepancy between model $q(\|y-\mu\|^2)$ and plug-in densities $q(\|y-\hat{\mu}\|^2)$ (i.e., $\int_{\mathbb{R}^d}  |q(\|y-\mu\|^2) \, - \, q(\|y-\hat{\mu}\|^2)| \, dy $),  with  $q$ unimodal and $Q$ the common c.d.f. of  the univariate components $Y_1, \ldots, Y_d$ with joint density $ q(\|y-\mu\|^2)$.
	
\noindent We proceed with a preparatory lemma which exploits the concavity of $\ell$, and which relates the difference in losses (\ref{loss3}) between estimates $\delta_g(X)=\delta_{0}(X) +(1-\omega)g(X)$ and $\delta_{0}(X)$, to the balanced squared-error loss difference.  Referring to the loss in (\ref{loss3}) as $L_{\omega,\ell}(\theta,\delta)$, we now define:
$$\Delta_{\omega,\ell}(\theta,\delta)=L_{\omega,\ell}(\theta,\delta)-L_{\omega,\ell} (\theta,\delta_{0}).$$

\begin{lem} (Lemma 6 in \cite{ms2020})
\label{l3.1}
Suppose that $X$ is spherically symmetric distributed about $\theta$ with density 
$f(\|x-\theta\|^2)$ . For the problem of estimating $\theta$ under (\ref{loss3}) with twice-differentiable, increasing, and concave $\ell$, we have
$$\Delta_{\omega,\ell}(\theta,\delta)\leq (1-\omega)^{2} \, \ell'\left\lbrace \,(1-\omega)\, \|\delta_{0}-\theta\|^2 \, \right\rbrace \, \Delta_{0,\ell}(\theta,\delta).$$
\end{lem}

	\noindent We now have the following.
\begin{theorem}
\label{theorem-ell-lebesgue}
Suppose that $X$ is spherically symmetric distributed about $\theta$ with density 
$f(\|x-\theta\|^2)$ and  that the function $\ell$ satisfies \textbf{C3}. For $d\geq 4$,  the estimator $\delta_{a,s}(X)$  in (\ref{Baranchik})  with condition (\ref{conditionsBaranchik}) dominates   $\delta_0(X)=X$ under loss  (\ref{loss3}) with $\delta_0(X)=X$ provided  $0<a < \displaystyle{2(d-2)}/\left({d \,E_0^{\star}\left[\displaystyle{\|Z\|^{-2}}\right]}\right)$,
where $$\displaystyle {Z \sim f^{\ast}(\|z-\ta \|^{2})=\displaystyle\frac{\ell'\{(1-\omega)\|z-\ta \|^{2}\} \, f(\|z-\ta \|^{2})}{\int_{\R^d}\ell'\{(1-\omega)\|z-\ta \|^{2}\} \, f(\|z-\ta \|^{2}) \, dz}}\,,$$	
\end{theorem}
and provided both $\mathbb{E}_0(\|X\|^2)$ and  $\mathbb{E}_0(\|X\|^{-2})$ are finite.

\noindent 	{\bf Proof.}   We show that the difference in risks between $\delta_0(X)=X$ and $\delta_{a,s}(X)$ 
 is non-negative under the given conditions. Let $\displaystyle K_Z=\int_{\R^d}\ell'\{(1-\omega)\|z-\ta \|^{2}\}f(\|z-\ta \|^{2})\, dz$, and $h^*$ be the radial density for $\|Z-\theta\|$.  Then, we have setting $g(X)\,=\, - a \frac{s(\|X\|^2)}{\|X\|^2} \, X\,$:
 \begin{eqnarray}
\label{eq-l-1}
	\Delta\cal{R}(\ta)
%	&=&	{\cal{R}}(\ta,\delta,\ell)-{\cal{R}}(\ta,\delta_{0},\ell)\nonumber\\
%	&=&E_{\ta}\left[ L_{\ta,\de_{0},\ell}(\theta,\de)-L_{\ta,\de_{0},\ell}(\theta,\de_{0})\right] \nonumber\\
	&=&
	E_{\ta}\left[ \Delta_{\ta,\ell}(\theta,\delta_{a,s})\right] \nonumber\\
	&\leq&
	E_{\ta}\left[ (1-\omega)^{2} \, \ell'\{(1-\omega) \, \|X-\theta\|^{2} \}\, \Delta_{0,\ell}(\theta,\delta_{a,s})\right] \nonumber\\
	&=&
	(1-\omega)^{2}E_{\ta}\left[ \ell'\{(1-\omega) \, \|X-\theta\|^{2} \} \, (\|\delta_{a,s}(X)-\ta\|^{2}-\|X-\ta\|^{2})\right] \nonumber\\
	&=&
\label{remarkbs1980}	(1-\omega)^{2}E_{\ta}\left[ \ell'\{(1-\omega) \, \|X-\theta\|^{2} \} \, \left( \|X\,+\,g(X)-\ta\|^{2}-\|X-\ta\|^{2}\right) \right] \\
	&=&
	\nonumber
	 K_Z \, (1-\omega)^{2} \, \left\lbrace E_{\ta}\left[ \, \|g(Z)\|^{2} \, + \, 2 (Z-\theta)^{\top} g(Z) \right] \right\rbrace
	 \\ & \leq & (1-\omega)^{2} \, \left\lbrace E_{\ta}\left[  \, \frac{a  s(\| Z\|^2)}{\| Z\|^2}\right] -2 \, \displaystyle\frac{d-2}{d}\int_{\R_{+}}r^{2}\displaystyle\int_{B_{r,\theta}}\frac{s(\| z\|^2)}{\|z\|^2} \, dV_{r,\theta}(z) \ h^*(r) \, dr \right\rbrace  \nonumber\\
	& = &
	a\, K_Z \, (1-\omega)^{2} E \left\lbrace\left( a \, \frac{1}{R^{2}}-2 \, \displaystyle\frac{d-2}{d}\right)
	E_{\ta}\left[R^2\frac{s(\| Z\|^2)}{\| Z\|^2}\Big{\arrowvert} \|Z-\ta\|=R\right]  \right\rbrace    \,,
	\end{eqnarray}
where  the first inequality follows from Lemma \ref{l3.1}, and the second inequality follows from a calculation using Lemma \ref{lem3.2} along with the defining inequalities $0 \leq s(\cdot) \leq 1$ and $s'(\cdot) \geq 0$. Finally, an application of Lemma \ref{lemBS1978} and the covariance inequality imply that 
	\begin{eqnarray}
	\Delta\cal{R}(\ta)
	&\leq & a \, K_Z \, (1-\omega)^{3}\left(a \, E_0 \left[{\|Z\|^{-2}}\right]  -2 \, \frac{d-2}{d}\right)E\left\{R^{2}E_{\ta}\left[\displaystyle\frac{s(\| Z\|^2)}{\| Z\|^2}\Big{\arrowvert} \|Z-\ta\|=R\right]\right\},   \nonumber
	\end{eqnarray}
	establishing the result.  \qed

\begin{rema}
A more direct proof of the above result is achieved by applying Theorem 2.1 of \cite{bs1980} immediately after (\ref{remarkbs1980}).  We have given here details following (\ref{remarkbs1980}) thus providing a more self-contained proof.  We also point out, as expected, that Theorems \ref{theorem-ell-lebesgue} and \ref{theorem-rho-lebesgue} match for $\ell(t)=\rho(t)=t$, and that Theorem \ref{theorem-ell-lebesgue} reduces to Theorem 2.1 of \cite{bs1980} for the unbalanced case $\omega=0$, as was the case for Theorem \ref{theorem-rho-lebesgue} in Section 2.

\end{rema}	

\subsection{Examples and discussion}
\label{subsection3.2}
As in Section 2, the dominance result above applies for many choices of $\ell$ and spherically symmetric densities $f$ guaranteeing the existence of Baranchik estimators that dominate the benchmark $\delta_0(X)=X$ under balanced loss $L_{\omega, \ell}$.  As in Section 2, Theorem \ref{theorem-ell-lebesgue}'s cut-off point on $a$ for scale mixture of normals is less than that of \cite{ms2020} (see subsection 3.2.1), but our dominance finding here is more generally applicable to all spherically symmetric densities subject to risk finiteness and to non completely monotone $\ell'$.

\noindent Theorem \ref{theorem-ell-lebesgue}'s cut-off point for dominance is representable in terms of 
$W=\|X-\theta\|^2$ as 
\begin{equation}
\label{b0}
a_0 \, = \, \frac{2(d-2)}{d} \, \frac{\mathbb{E}\left[ \, \ell'\{(1-\omega) W \}\right]} {\mathbb{E} \left[ \frac{\ell'\{(1-\omega)W) \}}{W}\right]} \,.
\end{equation}
An interesting case arises for $\ell(t)\,=\, t^q$ with $0 < q < 1$, with the above yielding 
\begin{equation}
\label{b0Lq}
a_0\,=\, a_0(q) \, = \,  \frac{2(d-2)}{d} \, \left\lbrace \mathbb{E}(W^{q-1})\right\rbrace/\left\lbrace \mathbb{E}(W^{q-2})\right\rbrace\,,
\end{equation}	
independently of $\omega$ as observed in \cite{ms2020}.  As expanded on in \cite{ms2020}, such a simple form of $a_0$ gives rise to simultaneous dominance with respect to both choices of $\ell$ and density $f$.  For instance, with 
$a_0(q)$ increasing for $q \in (0,1)$; which may be justified by writing $a_0(q)\,=\, \mathbb{E}_q(Z)$ with $Z$ having density to $z^{q-2} \, g_W(z)$ on $\mathbb{R}_+$ and observing that the family of such densities with parameter $q \in (0,1]$ has an increasing monotone likelihood ration in $Z$;  we have that $a_0 \geq a_0(q_0)$ for all $q \in [q_0,1]$, so that Theorem \ref{theorem-ell-lebesgue}'s Baranchik estimators $\delta_{a,s}(X)$ dominate $X$ for $a \leq a_0(q_0)$ simultaneously for all losses in (\ref{loss3}) with fixed $w \in [0,1)$ and such that $\ell(t)=t^q$ with $q \in [q_0,1]$.  We refer to \cite{ms2020} for such further examples with varying $f$ or $\ell$.

\subsubsection{Scale mixtures of normal distributions} 
	
\noindent  Theorem \ref{theorem-ell-lebesgue} applies for scale mixture of normals as in (\ref{mixtureofnormals}) for both completely monotone and non-completely monotone $\ell'$.  In the former case, Theorem \ref{theorem-ell-lebesgue}'s cut-off point, or equivalently (\ref{b0}), on $a$ for the Baranchik estimator $\delta_{a,s}(X)$ to dominate $X$ is smaller by a factor of $\frac{d-2}{d}$ that the one obtained in \cite{ms2020}, and it thus weaker.   For non-completely monotone $\ell'$, applications of Theorem \ref{theorem-rho-lebesgue} are novel however.  As an illustration, with $\ell(t) = 2 \Phi(t) -1$, the cut-off point  reduces to:
$$  a_0 \, = \, \frac{2(d-2)}{d} \,  \frac{\mathbb{E} \left[ e^{- \frac{(1-\omega)^2 \, W^2}{2}}\right]}{\mathbb{E} [W^{-1} \, e^{- \frac{(1-\omega)^2 \, W^2}{2}}]} \,,$$ 
for scale mixtures of normals, as well as for all spherically symmetric densities subject to the risk finiteness conditions.

\subsubsection{Example: Uniform distribution on a ball}

\noindent  Theorem \ref{theorem-ell-lebesgue} applies for the uniform distribution on the ball $B_{m,\theta}$ (i.e., $X \sim U(B_{m,\theta})$) with density in (\ref{density-uniform}).  Hence, the cut-off point $a_0$ can be evaluated using the density $ \frac{d}{2 m^d} \,w^{d/2-1} \, \mathbb{I}_{(0,m^2)}(w) $ for $W$. As an illustration for $\ell(t)\,=\, t^q$ with $0 < q < 1$, we obtain from (\ref{b0Lq}):

$$ a_0 \, = a_0(m) \, =\, \frac{2(d-2) \, m^2}{d} \, \; \frac{2q+d-4}{2q+d-2}\,.$$ 
Observe that the dominance finding has implications even in cases where $m$ is unknown,  but bounded below by a positive value $\underline{m}$, yielding that the dominance result $a \leq a_0(\underline{m})$, applicable for $X \sim U(B_{\underline{m},\theta})$, is robust to any discrepancy such that  $X \sim U(B_{m,\theta})$  with $m \geq \underline{m}$.  Moreover, the argument goes over when $m$ is random yielding a mixture of uniform distributions on balls as long as $\mathbb{P}(M \geq \underline{m})=1$, $M$ being the mixing parameter.

\subsubsection{Kotz type distribution}	

\noindent  As in Section 2, Theorem \ref{theorem-ell-lebesgue} applies to Kotz densities $f(\|x-\theta\|^2)$ with $f$ given in (\ref{model-kotz}).   Cut-off point $a_0$ in (\ref{b0}) is conveniently represented with $W=\|X-\theta\|^{2s} \sim \hbox{Gamma}(\nu,r)$.  For the case $\ell(t)=t^q$, a direct calculation of  (\ref{b0Lq}) yields for $\frac{q-2}{s}+\nu >0$
\begin{equation}
\nonumber
a_0 \, = \,  \frac{2(d-2)}{d}   \frac{\Gamma(\frac{q-1}{s}+\nu)}{\Gamma(\frac{q-2}{s}+\nu)}  \, r^{-1/s}\,.
\end{equation}

\section{Numerical illustration}

\subsection{Introduction}
We provide here frequentist risk evaluations illustrative of Theorem \ref{theorem-rho-lebesgue}, as well as provide details on the calculations themselves which we believe are useful for replication. The results of Section 3 can be illustrated in a similar fashion (but also see \cite{ms2020}).

\subsection{Calculation of risk under loss (\ref{loss2})}

The numerical evaluation of expectations with respect to a spherically symmetric density $f(\|x-\theta\|^2)$ may be expressed as a $d$ dimensional integral, but can be reduced to a two dimensional integral in our case of balanced risk function evaluations, and for estimators of $\theta$ that are equivariant under orthogonal transformations.
In this regard, the following lemma given in \cite{kariyaeaton1977} will turn out to be most useful.

\begin{lem} 
Let $X \sim f(\|x-\theta\|^2)$ with $X, \theta \in \mathbb{R}^d$, $\theta \neq 0$, $d>1$.    Let $\lambda=\|\theta\| $, $W=\|X\|^2$, and $T= \frac{\theta^{\top} X}{\lambda \|X\|}$.   Then, the joint density of $(T,W)$ is given by:

\begin{equation}
\label{pdftw}
\psi(t,w)\, = \, \frac{\pi^{\frac{d-1}{2}}}{\Gamma(\frac{d-1}{2})} \, w^{\frac{d}{2}-1} \, (1-t^2)^{\frac{d-3}{2}} \, f(w+\lambda^2-2 \lambda t w^{1/2})\,,
\end{equation}
for $t \in (-1,1)$ and $w>0$.  
\end{lem}
\noindent  Now, for equivariant estimators of $\theta$, which are of the form $h(\|X\|^2) X$ (e.g., \cite{eaton1989}), the loss in (\ref{loss2}) becomes
\begin{eqnarray*}
\omega \, \rho\left(\|h(\|X\|^2) X - X\|^2\right) \,  & + & \, (1- \omega) \, \rho\left(\|h(\|X\|^2) X - \theta\|^2\right) \\
\, =    \omega \, \rho\left( (h(W)-1)^2 \, W \right) \,  & + & \, (1- \omega) \, \rho\left(\lambda^2 + h^2(W) W \, - \, 2 \lambda \, T \, W^{1/2} \, h(W) \right) \,.
\end{eqnarray*}
Therefore, the associated frequentist risk for $\theta \in \mathbb{R}^d$  reduces to a function of $\lambda$ and can be evaluated as
\begin{equation}
\label{risknumerical}
 \omega \, \mathbb{E} \left\lbrace \rho\left( (h(W)-1)^2 \, W \right) \, \right\rbrace \,   +  \, (1- \omega) \, \mathbb{E} \left\lbrace \rho\left(\lambda^2 + h^2(W) W \, - \, 2 \lambda \, T \, W^{1/2} \, h(W) \right) \, \right\rbrace \, ,
\end{equation}
with the expectations taken with respect to density (\ref{pdftw}).
This illustrates the dimensional reduction to two dimensions and can be used to numerically evaluate the frequentist risk of Baranchik type estimators with $h$ of the form $h(w) \, = \, (1 - b \, \frac{s(w)}{w})$. 

\subsection{Illustration of Theorem \ref{theorem-rho-lebesgue}}

\noindent The dominance finding of Theorem \ref{theorem-rho-lebesgue} applies to Baranchik estimators $(1 - b \frac{s(\|X\|^2)}{\|X\|^2)}) X$ and expression (\ref{risknumerical}) can serve as a numerical illustration or comparison with inputs: {\bf (i)} $\omega \in [0,1)$ and $d \geq 4$; {\bf (ii)} a model density $f$; {\bf (iii)} a choice of $\rho$ satisfying condition {\bf C2}, taking $\rho'(0)=1$ without loss of generality; {\bf (iv)} a choice of $s(\cdot)$ satisfying (\ref{conditionsBaranchik}); and {\bf (v)} a choice of $b \leq
 a_0 (1-\omega)$ with $a_0$ given in (\ref{cutoff4}).

\noindent   We pursue by setting:  {\bf (i)} $\omega =1/2 $ and $d =6$; {\bf (ii)} a Kotz density as in (\ref{model-kotz}) with $r=s=1$, $\nu=4$ yielding the model density $X \sim f(\|x-\theta\|^2)$ with $f(t)\,=\, \frac{t e^{-t}}{3 \pi^3} \, \mathbb{I}_{(0,\infty)}(t)$; {\bf (iii)} $\rho(t) \, = \, \log(1+t)$; {\bf (iv)} a choice of $s(\cdot)$ satisfying (\ref{conditionsBaranchik}); and {\bf (v)} a choice of $b \leq a_0/2$  with expression (\ref{cutoff4}) yielding
\begin{equation}
\nonumber   \frac{a_0}{2} \, = \frac{1}{2}\,\frac{\frac{8}{3} \, I_3^2}{(\frac{2}{\pi^3} + I_3) I_2} \, \approx \, 0.595\,, 
\end{equation}
with $I_{k} = \frac{1}{3\pi^3}\int_0^{\infty} \frac{w^{k -1} e^{-w}}{1+w}  \, dw$, $I_3 \approx 0.01509$ and $I_2 \approx 0.00641$.

Figure 1 compares the risk of: {\bf (i)} $\delta_0(X)=X$, the Baranchik estimators {\bf (ii)} $s(t)=\frac{t}{1+t}$, $b=1/2$ (red in Figure 1), and {\bf (iii)} $s(t)=\frac{t}{1+t}$, $b=1$ (green in Figure 1), and the James-Stein estimator  {\bf (iv)}  with $s(t)=1$, and $b=0.5$ (blue in Figure 1).  The benchmark estimator $\delta_0$ is minimax (see first paragraph of Section 2)  with constant risk $R(\theta, \delta_0) = \frac{1}{2}\,\mathbb{E}_0(\log (1+\|X\|^2)) \,\approx \, 0.76606$, obtainable with $\|X-\theta\|^2 \sim \hbox{Gamma}(4,1)$ (see Subsection 2.3.4).   Both the first Baranchik and James-Stein estimators are minimax with sufficiently small cut-off point $b=0.5 \leq a_0/2$ as a consequence of Theorem \ref{theorem-rho-lebesgue}.  The maximal gains are attained at $\theta=0$ and are about $8.58\%$ and $10.43\%$.  The second Baranchik estimator has lower minimum risk, but its cut-off point is too large to satisfy the dominance condition of Theorem \ref{theorem-rho-lebesgue}.  Furthermore, the numerical evidence suggests that it is (barely) not minimax.  

 \begin{figure}[ht]
 \label{graphe1}
 \centering
\includegraphics[width=0.43\textwidth]{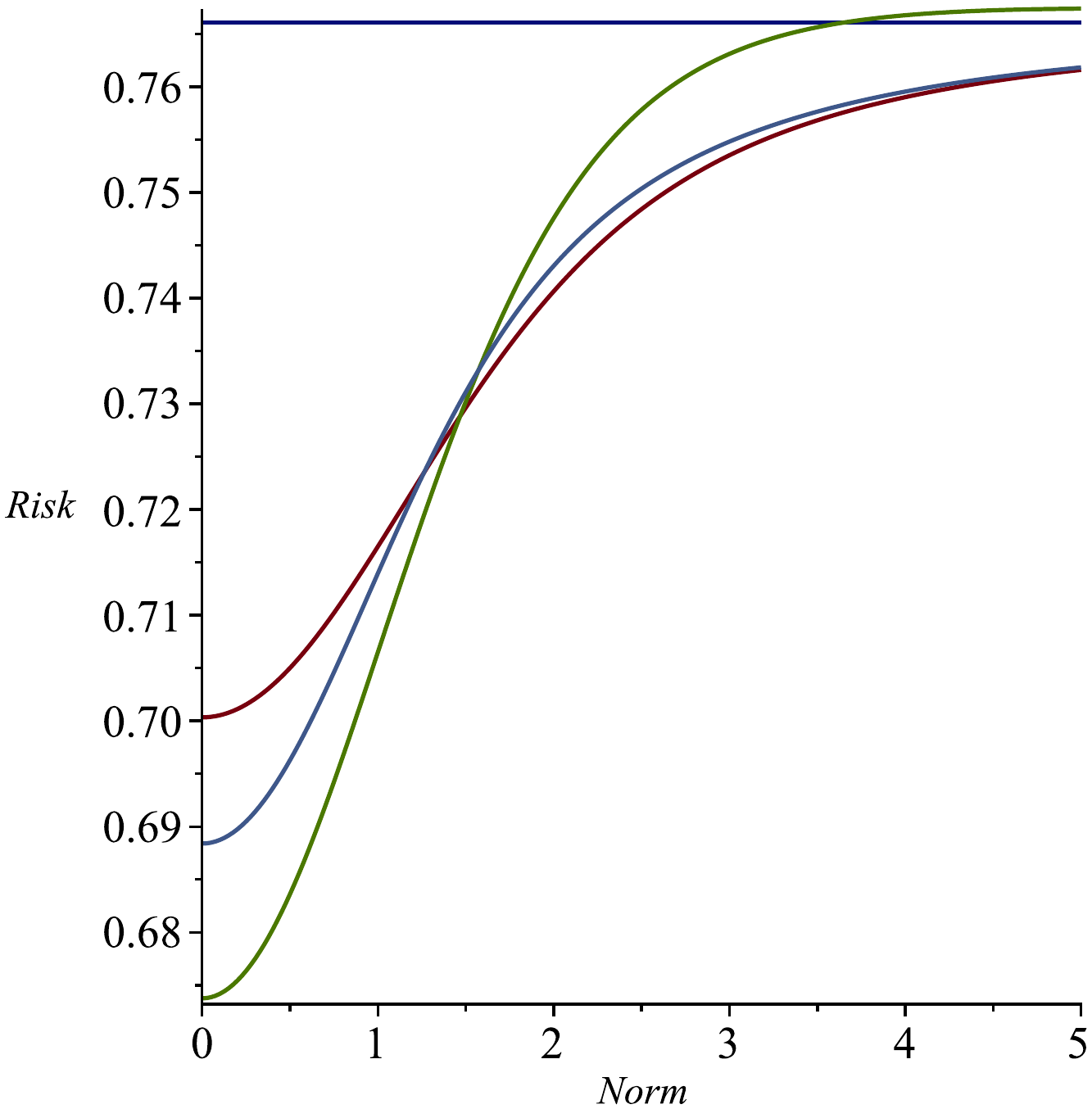}
  \caption{\footnotesize{  Frequentist risks as functions of $\|\theta\|$  of $\delta_0(X)=X$ (horizontal blue) and
  $\delta_{b,c}(X) \, = \, \left(1 - \frac{b}{\|X\|^2 + c}   \right) X$
  for $(b,c)=(0.5,1.0)$ (red), $(b,c)=(1.0,1.0)$ (green), and $(b,c)=(0.5,0.0)$ (blue); $X \sim \frac{1}{3 \pi^3}\,  \|x-\theta\|^2 \, e^{-\|x-\theta\|^2}$ (a Kotz density), balanced loss (\ref{loss2}) with $\rho(t) \, = \, \log(1+t)$, $\omega=1/2$.}} \label{graphe1}
 \end{figure}    
 	
\section{Concluding remarks}
For spherically symmetric distributed $X$ with densities $f(\|x-\theta\|^2)$, we  have provided frequentist risk improvements on the benchmark estimator $\delta_0(X)$ associated with balanced loss functions (\ref{loss2}) and (\ref{loss3}).  The findings, which apply to a wide class of Baranchik-type estimators, are unified with respect to the model $f$, with respect to the choices of $\rho$ in (\ref{loss2}) and $\ell$ in (\ref{loss3}), as well as the weight $\omega$ in these losses.    The findings extend earlier results of Marchand and Strawderman (\cite{ms2020}) which apply to scale mixtures of normals $f$ and completely monotone $\rho'$ or $\ell'$, while the results here apply to the whole class of spherically symmetric densities subject to risk finiteness, as well as to choices of $\rho$ and $\ell$ that are monotone increasing and concave.  

\noindent The findings testify to the ubiquitous nature of the effectiveness of shrinkage procedures as seen by the applicability of the dominance findings for a large class of densities $f$ and choices of $\rho, \ell$, and $\omega$.

\section*{Acknowledgements}
 
\'Eric Marchand's research is supported in part by the 
Natural Sciences and Engineering Research Council of Canada.

\section{Appendix}

Here are some technical results used in the paper. 

\begin{lem}
\label{lem A.1}
If $s: \mathbb{R}^d \to \mathbb{R}_+$ is a twice differentiable and concave function, then the function $x\rightarrow s(\|x\|^2)/\|x\|^2$ is superharmonic for $d\geq4$.
\end{lem}
\noindent {\bf{Proof.}}  Observe that $s(\cdot)$ must be non-decreasing.  A calculation of the Laplacian yields
\begin{equation}
\nonumber
\tri\left[s(\|x\|^2)/\|x\|^2\right] \,= \,  \frac{2}{\|x\|^4}\left[2\|x\|^4s''(\|x\|^2) +(d-4)\, \left(\|x\|^2s'(\|x\|^2)-s(\|x\|)^2\right)\right]\,.
\end{equation}
Since $s(0) \geq 0$ and $s$ is concave, we have  $s(0) \,+ ts'(t) \, \leq s(t)$ and the result follows.  \qed

The next result, referred to as a covariance inequality, is quite well known (e.g., Lemma 6.6, page 370 in \cite{lc1998}).
\begin{lem} \label{lemA.3} Let $Y$ be a random variable, and $g$ and $h$ be functions for which $\mathbb{E}[g(Y)]$, $\mathbb{E}[h(Y)]$, and $\mathbb{E}[g(Y) \, h(Y)]$ exist. 
\begin{enumerate}
\item[{\bf (a)}]  If one of the functions $g$ and $h$ is non-increasing and the other is non-decreasing, then $\mathbb{E}[g(Y) \, h(Y)] \leq 
\mathbb{E}[g(Y)] \,  \mathbb{E}[h(Y)]\,$;
\item[{\bf (b)}] 
If both $g$ and $h$ are either non-decreasing or non-increasing, then $\mathbb{E}[g(Y) \, h(Y)] \geq 
\mathbb{E}[g(Y)] \,  \mathbb{E}[h(Y)]\,$.
\end{enumerate}
\end{lem}
		
\noindent  The following is taken from \cite{fsw2018} and provides a useful decomposition for an expectation in terms of the radial distribution and uniform measures on balls.
	\begin{lem} 
	\label{lem3.2}
		Let $X$ have a spherically symmetric density about $\theta$, and  $g(X)$ be a weakly differentiable function such that $\mathbb{E}_{\theta}\left[ |(X-\theta)^{\top}g(X)|\right] <\infty$.  Then, assuming expectations exist, we have
		$$\mathbb{E}_{\theta}\left[(X-\theta)^{\top}g(X)\right] =\displaystyle\frac{1}{d} \, \mathbb{E} \, \left[ R^{2} \displaystyle\int_{B_{R,\theta}} div ( \, g(x) )\, dV_{R,\theta}(x)\right] $$
		where $\mathbb{E}$ denotes the expectation with respect to the radial distribution, and where $V_{r,\theta}(\cdot)$ is the 
		uniform distribution on $B_{R,\theta}$, the ball of radius $r$ centered at $\theta$.

	\end{lem} 
	{\bf{Proof.}}  See \cite{fsw2018}, Lemma 5.4, page 235.

		\begin{lem} (\cite{Plessis}, page 54)
			\label{lm3.3}
		Let $g$ be a superharmonic function, $Z_1$ have a uniform distribution on the  sphere $S_{r,\theta}$ centered  at $\theta$ with radius $r$, and $Z_2$ have a uniform distribution on the ball $B_{r,\theta}$ centered  at $\theta$ with radius $r$, then $\mathbb{E} \, g(Z_1) \leq \mathbb{E} \, g(Z_2) $.
		 
			% denoted $X\sim U(\|X-\theta\|= R^2)$, then
%			$$ E_{\theta}\left[g(X)\right] < E_{\theta}\left[g(Z)\right],$$ where $Z\sim U(\|Z-\theta\|^2\leq R^2$. That is 
%			$$ \frac{1}{A(S)}\int_S g(X) dA(X) \leq \frac{1}{M(B)}\int_Bg(Z) dM(Y),$$
%			where $A(S)$ and $M(B)$ represent the areas of the sphere $S$ and ball $B$, respectively.
		\end{lem}	
				
\begin{lem}
\label{lemBS1978}
Let $W$ be spherically symmetric Lebesgue density $g(\|w-\theta\|^2), \, w,\theta \in \mathbb{R}^d$.   Let $\beta: \mathbb{R}_+ \to \mathbb{R}_+$ be a non-decreasing function.  Then, $\mathbb{E}\left(  R^2 \, \frac{\beta(\|W\|^2)}{\|W\|^2} \big{|} \, \|W-\theta\|=R \right)$ is a non-decreasing function of $R>0$.
\end{lem}
{\bf Proof.}  See \cite{b+s1978}, pages 394-395 within the proof of their Theorem 3.3.1.  \qed

	\end{document}